\documentclass[10pt]{amsart}
\usepackage{amssymb}
\usepackage{amscd}

\def\today{\number\day\space\ifcase\month\or   January\or February\or
   March\or April\or May\or June\or   July\or August\or September\or
   October\or November\or December\fi\   \number\year}

\theoremstyle{definition}
\newtheorem{thm}{Theorem}[section]
\newtheorem{lem}[thm]{Lemma}
\newtheorem{prp}[thm]{Proposition}
\newtheorem{dfn}[thm]{Definition}
\newtheorem{cor}[thm]{Corollary}

\newtheorem{rmk}[thm]{Remark}

\newtheorem{exa}[thm]{Example}

\newcommand{\beq}{\begin{equation}}
\newcommand{\eeq}{\end{equation}}
\newcommand{\beqr}{\begin{eqnarray*}}
\newcommand{\eeqr}{\end{eqnarray*}}
\newcommand{\bal}{\begin{align*}}
\newcommand{\eal}{\end{align*}}
\newcommand{\bei}{\begin{itemize}}
\newcommand{\eei}{\end{itemize}}

\newcommand{\cP}{P}
\newcommand{\cK}{K}

\newcommand{\af}{\alpha}
\newcommand{\bt}{\beta}
\newcommand{\gm}{\gamma}
\newcommand{\dt}{\delta}
\newcommand{\ep}{\varepsilon}

\newcommand{\et}{\eta}

\newcommand{\io}{\iota}

\newcommand{\ld}{\lambda}
\newcommand{\sm}{\sigma}

\newcommand{\ph}{\varphi}

\newcommand{\rh}{\rho}

\newcommand{\Z}{{\mathbb{Z}}}

\newcommand{\C}{{\mathbb{C}}}
\newcommand{\N}{{\mathbb{N}}}

\newcommand{\tsr}{{\mathrm{tsr}}}
\newcommand{\Bsr}{{\mathrm{Bsr}}}

\newcommand{\andeqn}{\,\,\,\,\,\, {\mbox{and}} \,\,\,\,\,\,}

\newcommand{\ca}{C*-algebra}

\newcommand{\pj}{projection}

\newcommand{\ifo}{if and only if}

\newcommand{\mvnt}{Murray-von Neumann equivalent}
\newcommand{\mvnc}{Murray-von Neumann equivalence}
\newcommand{\nzp}{nonzero \pj}
\newcommand{\MIA}{M_{\infty} (A)}

\renewcommand{\S}{\subset}

\newcommand{\I}{\infty}

\newcommand{\TR}{{\mathrm{TR}}}
\newcommand{\Aut}{{\mathrm{Aut}}}
\newcommand{\card}{{\mathrm{card}}}
\newcommand{\Index}{{\mathrm{Index}}}
\newcommand{\id}{{\mathrm{id}}}
\newcommand{\End}{{\mathrm{End}}}
\newcommand{\Hom}{{\mathrm{Hom}}}

\title{Cancellation for inclusions of C*-algebras of finite depth}

\author{Ja A Jeong}

\address{Department of mathematics and
Research Institute of Mathematics,
  Seoul National University,
  Seoul, 151-747 Korea}

\email[]{jajeong@math.snu.ac.kr}

\author{Hiroyuki Osaka}

\address{Department of Mathematical Sciences,
  Ritsumeikan University, Kusatsu, Shiga,
  525-8577 Japan}

\email[]{osaka@se.ritsumei.ac.jp}

\author{N.~Christopher Phillips}

\address{Department of Mathematics, University of Oregon,
       Eugene OR 97403-1222, USA.}

\email[]{ncp@darkwing.uoregon.edu}

\author{Tamotsu Teruya}

\address{Department of Mathematical Sciences,
  Ritsumeikan University, Kusatsu, Shiga,
  525-8577 Japan}

\email[]{teruya@se.ritsumei.ac.jp}

\thanks{Research of the first author
partially supported by KRF-ABRL-R14-2003-006-01002-0.
Research of the second author
partially supported by the
Open Research Center Project for Private Universities: matching fund
{}from MEXT, 2004-2008, and by a Grant in Aid for Scientific Research,
Ritsumeikan University, 2005, 2006.
Research of the third author partially supported by
NSF grant DMS-0302401.}

\date{9~April 2007}

\begin{document}

\maketitle

\begin{abstract}
Let $1 \in A \subset B$ be a pair of \ca s with common
unit.
We prove that if $E \colon B \to A$ is a conditional
expectation with index-finite type and a quasi-basis of $n$
elements, then the topological stable rank satisfies
\[
\tsr (B) \leq \tsr (A) + n - 1.
\]
As an application we show that if an inclusion
$1 \in A \subset B$ of unital \ca s
has index-finite type and finite depth,
and $A$ is a simple
unital \ca\  with $\tsr (A) = 1$ and Property~(SP), then
$B$ has cancellation.
In particular, if $\af$ an action of
a finite group $G$ on $A,$ then the crossed product
$A \rtimes_{\af} G$ has cancellation.
For outer actions of $\Z,$
we obtain cancellation for $A \rtimes_{\af} {\mathbb{Z}}$
under the additional condition that $\af_* = \id$ on $K_0 (A).$
Examples are given.
\end{abstract}

\maketitle

\section{Introduction}\label{S:Intro}

For two projections $p,$ $q$ in a \ca, we write
$p \sim q$ if they are Murray-von Neumann equivalent.
A \ca~$A$ is said to have
{\emph{cancellation of projections}}
if whenever $p, q, r \in A$ are \pj s with $p \perp r,$ $q \perp r,$
and $p + r \sim q + r,$
then $p \sim q.$
If the matrix algebra $M_n (A)$ over $A$ has cancellation of
projections for each $n \in \N,$ we simply say that $A$ has
{\emph{cancellation}}.
Every \ca\  with cancellation is stably finite.

For a unital \ca~$A,$
if the topological stable rank $\tsr (A)$ of $A$
satisfies $\tsr (A) = 1,$
then $A$ has cancellation (Proposition 6.5.1 of~\cite{Bl1}).
For a stably finite simple \ca~$A,$
it has been a long standing open question,
settled negatively in~\cite{To},
whether cancellation implies $\tsr (A) = 1.$
The construction uses Villadsen's techniques~\cite{v}.
The example is also separable and nuclear.

Let $1 \in A \subset B$ be a unital inclusion
of \ca s with index-finite type and with finite depth.
In particular, $B$ could be a crossed product
$A \rtimes_{\af} G$ of a unital \ca\  by a finite group.
Our main result, Theorem~\ref{T:Mainindexfinite},
is that if $A$ is simple,
has topological stable rank~$1,$
and satisfies Property~(SP)
(every hereditary C*-subalgebra contains a nonzero projection),
then $B$ has cancellation.

As a corollary, suppose that
$A$ is a simple unital \ca\  with tracial topological rank zero
($\TR (A) = 0$; Definition 3.6.2 of~\cite{Hl}),
and $\af \colon G \to \Aut (A)$ is an action of a
finite group~$G$ on $A.$
Then $A \rtimes_{\af} G$ has cancellation.
Examples for $A$ include all
simple unital AH-algebras with real rank zero
and slow dimension growth.
Recently, the third author (\cite{Ph11}, \cite{PhtRp1a})
has proved that if in addition $\af$ has the tracial Rokhlin property,
then $\TR (A \rtimes_{\af} G) = 0.$
It follows that, in this case, $\tsr (A \rtimes_{\af} G) = 1.$
The result of this paper shows that no conditions on the action
are needed for cancellation.

As an intermediate result, we prove in Section~\ref{S:tsr}
that if $1 \in A \subset B$
is a unital inclusion of \ca s,
and if there is a faithful conditional
expectation $E \colon B \to A$ with index-finite type and a
quasi-basis of $n$ elements (\cite{wata}; detailed reference below),
then
\[
\tsr (B) \leq \tsr (A) + n - 1
\]
In particular, this applies if $B = A \rtimes_{\af} G$
and $\card (G) = n.$

Another important ingredient is a result of
Blackadar (Theorem~\ref{T:BCanc} below):
if $B$ is a simple \ca,
and $\cP$ is a set of nonzero \pj s in $B$
containing, in a suitable sense, arbitrarily small \pj s,
such that $\sup_{p\in \cP} \tsr (p B p) < \infty,$
then $B$ has cancellation.

The authors would like to thank Ken Goodearl for his
kind permission to present his proof of
Theorem~\ref{T:BCanc}.
The second author would like to thank Ken Goodearl and Masaru Nagisa
for fruitful discussions.

\section{Topological stable rank}\label{S:tsr}

For a unital \ca~$A,$ recall
that the {\emph{topological stable rank}} $\tsr (A)$ of $A$ is defined
to be the least integer $n$
such that the set ${\mathrm{Lg}}_n (A)$
of all $n$-tuples $(a_1,, a_2 \dots, a_n) \in A^n$
which generate $A$ as a left ideal is dense in $A^n.$
(See Definition~1.4 of~\cite{Rf1}.)
The topological stable rank of a nonunital
\ca\  is defined to be that of its smallest unitization.
Note that $\tsr (A) = 1$ is equivalent to
density of the set of invertible elements in $A.$
Furthermore, $\tsr (A) = 1$ implies that $\tsr (A \otimes M_{n}) = 1$
for all $n,$ and that $\tsr (A \otimes \cK) = 1,$
where $\cK$ is the algebra of compact
operators on a separable infinite dimensional Hilbert space.
Conversely, if $\tsr (A \otimes M_{n}) = 1$ for some $n,$
or if $\tsr (A \otimes \cK) = 1,$
then $\tsr (A) = 1.$
(See Theorems 3.3 and~3.6 of~\cite{Rf1}.)
Simple AH~algebras with slow dimension growth have
topological stable rank one (Theorem~1 of~\cite{BDR}),
as do irrational rotation algebras (\cite{pu}).
If $A$ is unital and $\tsr (A) = 1,$
the $A$ has cancellation (Proposition 6.5.1 of~\cite{Bl1}).
It follows immediately that $A$ is stably finite in the sense
that no matrix algebra $M_n (A)$ has an infinite projection.

As in~\cite{wata}
(see Definition 1.2.1 and the following discussion there),
if $1 \in A \subset B$ is a pair of \ca s with a common unit,
then a {\emph{conditional expectation}}
$E \colon B \to A$ is a
positive $A$-bimodule map of norm one.
Following Definition~1.2.2 and Lemma 2.1.6 of~\cite{wata},
if $E$ is faithful (a standing assumption in~\cite{wata};
see the discussion after Definition 1.2.1 there),
a {\emph{quasi-basis}} for $E$ is a finite family
$\big( (u_1, u_1^*), \, (u_2, u_2^*), \, \dots, \, (u_n, u_n^*) \big)$
of elements of $B \times B$
such that
\[
b = \sum_{j = 1}^n u_j E (u_j^* b)
  = \sum_{j = 1}^n E (b u_j) u_j^*
\]
for all $b \in B,$
the expectation $E$ has {\emph{index-finite type}} if $E$
has a quasi-basis,
and the index of $E$ is then defined by
$\Index (E) = \sum_{j = 1}^n u_j u_j^*.$
By Proposition~1.2.8 and Lemma~2.3.1 of~\cite{wata},
the index is a positive invertible
central element of $B$ that does not depend on the choice of the
quasi-basis.
In particular, if $1 \in A \subset B$ is a pair of simple unital \ca s,
then $\Index (E)$ is a positive scalar.
By abuse of language, we will say that $1 \in A \subset B$
has index-finite type if there is a faithful
conditional expectation $E \colon B \to A$ with index-finite type.

The following example is standard.
(See Lemma~3.1 of~\cite{OT} for a good deal more.)

\begin{exa}\label{E:InclCrPrd}
Let $A$ be a unital \ca,
let $G$ be a finite group,
and let $\af \colon G \to \Aut (A)$ be an action of $G$ on $A.$
For $g \in G,$
let $u_g \in A \rtimes_{\af} G$ be the standard unitary
in the crossed product, implementing $\af_g.$
Then the function $E \colon A \rtimes_{\af} G \to A,$
given by $E \big( \sum_{g \in G} a_g u_g \big) = a_1,$
is a conditional expectation with index-finite type,
$\big( (u_{g}, u_{g}^{*}) \big)_{g \in G}$ is a quasi-basis for $E,$
and $\Index (E) = \card(G) \cdot 1_{A \rtimes_{\af} G}.$
\end{exa}

\begin{thm}\label{T:Main}
Let $B$ be a unital \ca,
let $A \S B$ be a unital subalgebra,
let $E \colon B \to A$ be
a faithful conditional expectation with index-finite type,
and let $\big( (v_k, v_k^*) \big)_{1 \leq k \leq n}$
be a quasi-basis for $E.$
Then $\tsr (B) \leq \tsr (A) + n - 1.$
\end{thm}

\begin{proof}
Set $m = \tsr (A) - 1.$
We prove that ${\mathrm{Lg}}_{m + n} (B)$ is dense in $B^{m + n}.$
Let $b_1, b_2, \ldots, b_{m + n} \in B,$
and let $\ep > 0.$
Write
\[
b_j = \sum_{k = 1}^n a_{j, k} v_k
\]
for $1 \leq j \leq m + n,$ with all $a_{j, k} \in A.$

We will work with various sizes of nonsquare matrices over $A$
and over $B.$
We regard all of them as elements of the \ca\  $B \otimes {\cK}$
by placing each one in the upper left corner of an infinite
matrix, and taking all remaining entries of the infinite
matrix to be zero.
Multiplication of matrices of compatible sizes thus
becomes the usual multiplication in $B \otimes {\cK},$
and all the usual properties of multiplication in \ca s are valid.
We also write $1_l$ for the $l \times l$ identity matrix,
which, according to the convention just made,
is a \pj\  in $B \otimes {\cK}.$

Set
\[
a = \left( \begin{array}{ccccccc}
  a_{1, 1}         &  a_{1, 2}         &  \hdots  &  a_{1, n}      \\
  a_{2, 1}         &  a_{2, 2}         &  \hdots  &  a_{2, 2}      \\
  \vdots           &  \vdots           &          &  \vdots        \\
  a_{m + n, \, 1}  &  a_{m + n, \, 2}  &  \hdots  &  a_{m + n, \, n}
\end{array} \right)
\andeqn
v = \left( \begin{array}{c}
  v_1     \\
  v_2     \\
  \vdots  \\
  v_n
\end{array} \right).
\]
This gives
\[
a v = \left( \begin{array}{c}
  b_1     \\
  b_2     \\
  \vdots  \\
  b_{m + n}
\end{array} \right).
\]

According to Definition~6.2 and Lemma~6.3 of~\cite{Rf1},
the \ca~$A$ satisfies the property $L_m (n),$
so that there exists
\[
x = \left( \begin{array}{ccccccc}
  x_{1, 1}         &  x_{1, 2}         &  \hdots  &  x_{1, n}      \\
  x_{2, 1}         &  x_{2, 2}         &  \hdots  &  x_{2, 2}      \\
  \vdots           &  \vdots           &          &  \vdots        \\
  x_{m + n, \, 1}  &  x_{m + n, \, 2}  &  \hdots  &  x_{m + n, \, n}
\end{array} \right),
\]
with $x_{j, k} \in A$ for
$1 \leq j \leq m + n$ and $1 \leq k \leq n,$
such that $\| x - a \| < \ep/\|v\|$ and
$x$ is left invertible.
This last condition means that there is a $n \times (m + n)$ matrix $z,$
with entries in $A,$
such that $z x = 1_n.$
Then
$1_n = x^* z^* z x \leq \| z^* z \| x^* x.$
Thus, with $\dt = \| z^* z \|^{ - 1},$
we have $x^* x \geq \dt \cdot 1_n.$

Define
\[
y_j = \sum_{k = 1}^n x_{j, k} v_k \in B
\]
for $1 \leq j \leq m + n.$
Then
\[
x v = \left( \begin{array}{c}
  y_1     \\
  y_2     \\
  \vdots  \\
  y_{m + n}
\end{array} \right).
\]
We therefore get the following relation,
which by our convention
is really an inequality between matrices in $B \otimes {\cK}$
whose only nonzero entries are in the $1, 1$ position:
\[
\sum_{j = 1}^{m + n} y_j^* y_j
 = v^* x^* x v
\geq \dt v^* v
 = \dt \sum_{k = 1}^n v_k^* v_k
 = \dt \cdot {\mathrm{Index}} (E).
\]
The inequality is still correct when interpreted in $B.$
Since ${\mathrm{Index}} (E)$ is a positive invertible element of $B,$
it follows that $\sum_{j = 1}^{m + n} y_j^* y_j$
is invertible in $B.$
Therefore $(y_1, y_2, \ldots, y_{m + n}) \in {\mathrm{Lg}}_{m + n} (B).$
We have $\| x v - a v \| \leq \| x - a \|\|v\| < \ep,$
so that $\| y_k - b_k \| < \ep$ for $1 \leq j \leq m + n.$
This proves that ${\mathrm{Lg}}_{m + n} (B)$ is dense in $B^{m + n}.$
\end{proof}

Using Theorem~\ref{T:Main} we can sharpen the estimate
of Proposition~5.3 of~\cite{OT}.

\begin{cor}\label{Cor:Proposition5.3}
Let $B$ be a unital \ca,
let $A \S B$ be a unital subalgebra,
let $E \colon B \to A$ be
a faithful conditional expectation with index-finite type,
and let $\big( (v_k, v_k^*) \big)_{1 \leq k \leq n}$
be a quasi-basis for $E.$
Then
\[
\tsr (A) \leq n \cdot \tsr (B) + n^2 - 2 n + 1.
\]
\end{cor}

\begin{proof}
In the proof of Proposition~5.3 of~\cite{OT},
substitute the estimate of Theorem~\ref{T:Main}
for the estimate of Corollary~2.6 of~\cite{OT}.
\end{proof}

\begin{thm}\label{T:Gp}
Let $A$ be a \ca, and let $\af \colon G \to \Aut (A)$
be an action of a finite group $G$ on $A.$
Then $\tsr (A \rtimes_{\af} G) \leq \tsr (A) + \card (G) - 1.$
\end{thm}

\begin{proof}
The unital case follows from Theorem~\ref{T:Main}
and Example~\ref{E:InclCrPrd}.
For the nonunital case,
let $A^+$ be the unitization of $A,$ and
observe that $A \rtimes_{\af} G$ is an ideal in $A^+ \rtimes_{\af} G.$
Using, in order, Theorem~4.4 of~\cite{Rf1},
the unital case, and Definition~1.4 of~\cite{Rf1}, we get
\[
\tsr (A \rtimes_{\af} G)
\leq \tsr (A^+ \rtimes_{\af} G)
\leq \tsr (A^+) + \card (G) - 1
 = \tsr (A) + \card (G) - 1,
\]
as desired.
\end{proof}

\begin{rmk}\label{R:CycGp}
As pointed out in Example 8.2.1 of~\cite{bl3},
Theorems~4.3 and~7.1 of~\cite{Rf1} can be used to show that
for any action $\af \colon {\mathbb{Z}} / 2 {\mathbb{Z}} \to \Aut (A),$
one has
$\tsr ( A \rtimes_{\af} {\mathbb{Z}} / 2 {\mathbb{Z}} )
 \leq \tsr (A) + 1.$
The point is that $A \rtimes_{\af} {\mathbb{Z}} / 2 {\mathbb{Z}}$
is a quotient of $A \rtimes_{\af} {\mathbb{Z}}.$
The argument actually works for any finite cyclic group.
This estimate is the same as that of Theorem~\ref{T:Gp}
for ${\mathbb{Z}} / 2 {\mathbb{Z}},$
and better for other cyclic groups.
\end{rmk}

\begin{rmk}
The estimate in Theorem~\ref{T:Gp} is the best possible
of its form.
There is a (nonsimple) unital \ca~$A$
with $\tsr (A) = 1$ and an action
$\af \colon {\mathbb{Z}} / 2 {\mathbb{Z}} \to \Aut (A)$
such that $\tsr (A \rtimes_{\af} {\mathbb{Z}} / 2 {\mathbb{Z}}) = 2.$
See Example 8.2.1 of~\cite{bl3}.
\end{rmk}

\section{Inclusions of C*-algebras of finite depth}

The notion of finite depth for subfactors is well known.
(See, for example, Section~4.6 of~\cite{GHJ}.)
The basic properties of inclusions of \ca s with finite depth
are similar, but have not appeared in the literature.
They do not differ greatly from the subfactor case,
but, in the interest of completeness, we give proofs here.

\begin{dfn}\label{D:FiniteDepth}
Let $1 \in A \subset B$ be an inclusion of unital \ca s
with a conditional expectation $E \colon B \to A$
with index-finite type.
Set $B_0 = A,$ $B_1 = B,$
and $E_1 = E.$
Recall the \ca\  version of the basic construction
(Definition~2.2.10 of~\cite{wata}, where it is called the
C*~basic construction).
We inductively define $e_k = e_{B_{k - 1}}$
and $B_{k + 1} = C^* (B_k, e_k),$
the Jones projection and \ca\  for the basic construction
applied to $E_k \colon B_k \to B_{k - 1},$
and take $E_{k + 1} \colon B_{k + 1} \to B_k$ to be
the dual conditional expectation $E_{B_k}$ of
Definition~2.3.3 of~\cite{wata}.
This gives the {\emph{tower of iterated basic constructions}}
\[
B_0 \subset B_1 \subset B_2 \subset \cdots \subset B_k \subset \cdots,
\]
with $B_0 = A$ and $B_1 = B.$
It follows from Proposition 2.10.11 of~\cite{wata} that
this tower does not depend on the choice of $E.$

We then say that the inclusion $A \subset B$
has {\emph{finite depth}}
if there is
$n \in \N$ such that
$(A' \cap B_{n}) e_{n} (A' \cap B_{n}) = A' \cap B_{n + 1}.$
We call the least such $n$ the {\emph{depth}} of the inclusion.
\end{dfn}

\begin{dfn}\label{D:MinExpt}
Let $1 \in A \subset B$ be an inclusion of unital \ca s
with a conditional expectation $E \colon B \to A$
with index-finite type.
We say that $E$ is {\emph{pseudominimal}} if
$\Index (E)$ is a scalar multiple of $1_A$
and $E (c b) = E (b c)$
for all $c \in A' \cap B$ and $b \in B.$
\end{dfn}

When $A$ and $B$ have trivial centers,
Lemma~3.11 of~\cite{kw} shows that
a minimal conditional expectation is pseudominimal.

The following is an extended version of Lemma~3.11 of~\cite{kw}.

\begin{lem}\label{MinimalConditional}
Let $1 \in A \subset B$ be an inclusion of C*-algebras,
with conditional expectation $E \colon B \to A$
with index-finite type.
Suppose that $A$ is simple.
Then there exists a pseudominimal conditional expectation
$F \colon B \rightarrow A$ with index-finite type.
\end{lem}

\begin{proof}
Using Theorem~3.3 of~\cite{Iz}
(also see Lemma~2.2 and Remark~2.4(3) of~\cite{Iz}),
there exist central projections $p_1, p_2, \ldots, p_n \in B$
such that $B p_j$ is simple for $1 \leq j \leq n$
and $B = \bigoplus B p_j.$
We then have unital inclusions
$p_j \in p_j A p_j \subset p_j B p_j = B p_j$ for each $j.$
Note that $p_j B p_j = B p_j$ is a simple ideal in $B.$
However, $p_j A p_j$ is not a subalgebra of $A,$
only of $B.$
In fact, since $A$ is simple,
$\ph_j (a) = a p_j$ defines an isomorphism from $A$ to $p_j A p_j.$

There are conditional expectations $E_j \colon p_j B p_j \to p_j A p_j$
given by $E_j (b) = p_j E (b) p_j$ for $b \in p_j B p_j.$
If
$\big( (u_k, u_k^*) \big)_{1 \leq k \leq m}$
is a quasi-basis for $E,$
then
$\big( (p_j u_k, \, (p_j u_k)^*) \big)_{1 \leq k \leq m}$
is a quasi-basis for $E_j.$
Thus, the unital inclusion $p_j \in p_j A p_j \subset p_j B p_j$
has index-finite type.
Following Proposition~2 of~\cite{KwW},
let $F_j \colon p_j B p_j \to p_j A p_j$
be the minimal conditional expectation from $p_j B p_j$
onto $p_j A p_j$
(which automatically also has index-finite type).
Its index can be written
as $\Index (F_j) = \ld_j p_j$ because it is in the center of
the simple \ca\  $p_j B p_j.$
Let
$\big( (v_{j, l}, v_{j, l}^*) \big)_{1 \leq l \leq m_j}$
be a quasi-basis for~$F_j.$
Choose $\bt_1, \bt_2, \ldots, \bt_n > 0$
such that
\[
\sum_{j = 1}^n \bt_j = 1
\andeqn
\bt_1^{-1} \ld_1
 = \bt_2^{-1} \ld_2 = \cdots = \bt_n^{-1} \ld_n.
\]
Define
$F \colon B \to A$ by
\[
F (b) = \sum_{j = 1}^n \bt_j \ph_j^{-1} ( F_j (p_j b)).
\]
Then $F$ is a conditional expectation.

Set $w_{j, l} = \bt_j^{-1/2} v_{j, l}.$
We claim that
$\big( ( w_{j, l}, w_{j, l}^*) \big)_{1 \leq j \leq n,
                                 \, 1 \leq l \leq m_j}$
is a quasi-basis for~$F.$
First note that for $a \in p_j A p_j$ we have
$v_{j, l} \ph_j^{-1} (a) = v_{j, l} a$
and $\ph_j^{-1} (a) v_{j, l} = a v_{j, l},$
since $v_{j, l} \in p_j B p_j.$
Then for $b \in B$ we have,
using $p_k v_{j, l} = \dt_{j, k} v_{j, l}$ at the second step,
\begin{align*}
\sum_{j = 1}^n \sum_{l = 1}^{m_j} w_{j, l} F (w_{j, l}^* b)
& = \sum_{j = 1}^n \sum_{l = 1}^{m_j} \bt_{j}^{-1}
      v_{j, l} \sum_{k = 1}^n \bt_k \ph_k^{-1} ( F_k (p_k v_{j, l}^* b))
       \\
& = \sum_{j = 1}^n \sum_{l = 1}^{m_j}
      v_{j, l} \ph_j^{-1} ( F_j (v_{j, l}^* b p_j))
       \\
& = \sum_{j = 1}^n \sum_{l = 1}^{m_j}
      v_{j, l} F_j (v_{j, l}^* b p_j)
  = \sum_{j = 1}^n b p_j
  = b.
\end{align*}
The proof that
$\sum_{j = 1}^n \sum_{l = 1}^{m_j} F (b w_{j, l}) w_{j, l}^* = b$
is similar.

Now we check the conditions in the definition of pseudominimality.
By construction, 
we have
\[
\Index (F)
  = \sum_{j = 1}^n \sum_{l = 1}^{m_j} w_{j, l} w_{j, l}^*
  = \sum_{j = 1}^n \bt_j^{-1} \ld_j p_j,
\]
which is a scalar by the choice of the $\bt_j.$
For the commutation relation, let $c \in A' \cap B$ and $b \in B.$
Then $p_j c \in (p_j A p_j)' \cap (p_j B p_j),$ so
minimality of $F_j,$ Lemma~3.11 of~\cite{kw}, and
centrality of $p_j,$ imply $F_j (p_j b c) = F_j (p_j c b).$
Thus
\[
F (b c) = \sum_{j = 1}^n \bt_j \ph_j^{-1} ( F_j (p_j b c))
  = \sum_{j = 1}^n \bt_j \ph_j^{-1} ( F_j (p_j c b))
  = F (c b).
\]
This completes the proof.
\end{proof}

The next proposition was proved in
the II$_1$ factor case by Pimsner and Popa in
Theorem~2.6 of~\cite{PP}.

\begin{prp}\label{tower}
Let $1 \in A \subset B$ be an inclusion of unital \ca s
with a faithful conditional expectation $E \colon B \to A$
with index-finite type
such that $\Index (E) \in A.$
Using the notation of Definition~\ref{D:FiniteDepth},
for $n \geq 1$
we can identify $B_{2 n}$ as the basic construction for
$A \subset B_n$ and the conditional expectation
$F = E_1 \circ E_2 \circ \cdots \circ E_n.$
\end{prp}

\begin{proof}
As in the notation of Definition~\ref{D:FiniteDepth},
let $e_k$ be the Jones \pj\  for the inclusion of $B_{k - 1}$ in $B_k,$
so that $B_{k + 1} = C^* (B_k, e_k).$
Further, set $z = \Index (E).$
By Lemma~2.3.1 and Proposition~2.3.4 of~\cite{wata},
and induction,
$z$ is a positive invertible element which is in the center of $B_n$
for all $n \geq 0.$

Following the beginning of Section~2 of~\cite{PP}
(but noting that our indexing conventions differ),
define $g_n^{(k)} \in B_k' \cap B_{2 n + k}$ by
\[
g_n^{(k)}
 = (e_{n + k} e_{n + k - 1} \cdots e_{k + 1})
      (e_{n + k + 1} e_{n + k} \cdots e_{k + 2})
      \cdots (e_{2 n + k - 1} e_{2 n + k - 2} \cdots e_{n + k}).
\]
(Take $g_0^{(k)} = 1.$)
Then set
$f_n^{(k)} = z^{n (n - 1)/2} g_n^{(k)}.$
Using, say, Definition~2.3.3,
Propositions~2.1.1, 2.3.2, and~2.3.4, and Lemma~2.3.5 of~\cite{wata},
we have the usual relations among the \pj s $e_k$:
\begin{enumerate}
\item\label{Jones1}
$e_k x e_k = E_k (x) e_k$ for $x \in B_k.$
\item\label{Jones2}
$e_k e_{k \pm 1} e_k = z^{-1} e_k.$
\item\label{Jones3}
$e_k e_l = e_l e_k$ when $| k - l | > 1.$
\end{enumerate}
Since $z$ is a positive, invertible, and central,
the arguments of Section~2 of~\cite{PP}
(through the calculations in the proof of Theorem~2.6 there)
go through (replacing scalars by powers of $z$),
and imply that $f_n = f_n^{(0)}$ is a \pj,
and (using our indexing) that
$E_{2 n + k} \big( g_n^{(k)} \big)
   = z^{-n} g_{n - 1}^{(k + 1)}$
for all $k$ and $n.$
It follows that
$E_{n + 1} \circ E_{n + 2} \circ \cdots \circ E_{2 n} (f_n) = z^{-n}.$

We now want to calculate $f_n x f_n$ for $x \in B_n.$
Since $f_n^* = f_n,$ we start out as follows,
using $e_n x e_n = E_n (x) e_n$ at the first step and
\[
E_n (x) \in B_{n - 1}
 \S B_{2 n}
    \cap \{ e_{n + 1}, \, \ldots, \, e_{2 n - 2} , \, e_{2 n - 1} \}'
\]
at the second step:
\begin{align*}
& (e_{2 n - 1} e_{2 n - 2} \cdots e_{n})
    x (e_{n} \cdots e_{2 n - 2} e_{2 n - 1})
    \\
& \hspace*{3em} {\mbox{}}
    = (e_{2 n - 1} e_{2 n - 2} \cdots e_{n + 2})
     E_n (x) e_{n + 1} (e_{n + 2} \cdots e_{2 n - 2} e_{2 n - 1})
    \\
& \hspace*{3em} {\mbox{}}
    = E_n (x) (e_{2 n - 1} e_{2 n - 2} \cdots e_{n + 3})
        (e_{n + 2} e_{n + 1} e_{n + 2})
        (e_{n + 3} \cdots e_{2 n - 2} e_{2 n - 1})
    \\
& \hspace*{3em} {\mbox{}}
    = E_n (x) (e_{2 n - 1} e_{2 n - 2} \cdots e_{n + 3})
        z^{- 1} e_{n + 2}
        (e_{n + 3} \cdots e_{2 n - 2} e_{2 n - 1}).
\end{align*}
Iterating the last two steps gives
\[
(e_{2 n - 1} e_{2 n - 2} \cdots e_{n})
    x (e_{n} \cdots e_{2 n - 2} e_{2 n - 1})
  = z^{- (n - 1)} E_n (x) e_{2 n - 1}.
\]
Set $x_{n - 1} = z^{- (n - 1)} E_n (x) \in B_{n - 1}.$
We now calculate
\[
(e_{2 n - 2} e_{2 n - 3} \cdots e_{n - 1})
       x_{n - 1} e_{2 n - 1} (e_{n - 1} \cdots e_{2 n - 3} e_{2 n - 2})
\]
by a similar method.
First, use the commutation relation~(\ref{Jones3}) to write
this expression as
\[
e_{2 n - 2} \big[ (e_{2 n - 3} e_{2 n - 4} \cdots e_{n - 1})
     x_{n - 1} (e_{n - 1} \cdots e_{2 n - 4} e_{2 n - 3}) \big]
               e_{2 n - 1} e_{2 n - 2}.
\]
The method used above shows that the term in brackets is equal to
\[
z^{- (n - 2)} E_{n - 1} (x_{n - 1}) e_{2 n - 3},
\]
from which it follows, using the fact that $e_{2 n - 2}$
commutes with $B_{n - 2}$ at the second step, that
\begin{align*}
& (e_{2 n - 2} e_{2 n - 3} \cdots e_{n - 1})
       x_{n - 1} e_{2 n - 1} (e_{n - 1} \cdots e_{2 n - 3} e_{2 n - 2})
    \\
& \hspace*{3em} {\mbox{}}
   = e_{2 n - 2}
  \big[ z^{- (n - 2)} E_{n - 1} (x_{n - 1}) e_{2 n - 3} \big]
               e_{2 n - 1} e_{2 n - 2}
    \\
& \hspace*{3em} {\mbox{}}
    = z^{- (n - 2)} E_{n - 1} (x_{n - 1})
           (e_{2 n - 2} e_{2 n - 3}) ( e_{2 n - 1} e_{2 n - 2} ).
\end{align*}
Proceeding inductively, and putting in the appropriate power of $z,$
we finally arrive at
\[
f_n x f_n
 = \left(E_1 \circ E_2 \circ \cdots \circ E_n \right)(x) f_n.
\]
This result implies in particular that the map $x \mapsto x f_n$
is injective on $B_0.$
Proposition~2.2.11 of~\cite{wata} therefore implies that the
subalgebra $C^* (B_n, f_n) \S B_{2 n}$ is the basic construction for
$B_0 \subset B_n.$

Using $f_n^{(k)}$ in place of $f_n = f_n^{(0)},$
we obtain analogous results,
and, in particular, for every $k,$ the
subalgebra $C^* \big( B_{n + k}, \, f_n^{(k)} \big) \S B_{2 n + k}$
is the basic construction for $B_k \subset B_{n + k}.$

We now prove by induction on $n$ that
$C^* \big( B_{n + k}, \, f_n^{(k)} \big) = B_{2 n + k}$ for all $k.$
This is true by hypothesis for $n = 1,$
so suppose it is known for $n - 1.$
The inclusion $B_{k} \S B_{k + 1},$
with conditional expectation $E_{k + 1},$
has index-finite type by Proposition 1.6.6
and Definition~2.3.3 of~\cite{wata}.
Let
$\big( (v_1, v_1^*), \, (v_2, v_2^*), \, \dots, \, (v_r, v_r^*) \big)$
be a quasi-basis for this inclusion.
Then
\[
\sum_{j = 1}^r v_j e_{k + 1} v_j^*
 = \sum_{j = 1}^r E_{k + 1} (v_j) v_j^* = 1,
\]
and $v_j$ commutes with $e_l$ for $l > k + 1,$ so
\begin{align*}
& \sum_{j = 1}^r v_j g_n^{(k)} v_j^*        \\
& \hspace*{1em} {\mbox{}}
 = (e_{n + k} e_{n + k - 1} \cdots e_{k + 2})
      (e_{n + k + 1} e_{n + k} \cdots e_{k + 2})
      \cdots (e_{2 n + k - 1} e_{2 n + k - 2} \cdots e_{n + k}).
\end{align*}
(The difference from the definition of $g_n^{(k)}$ is that here the
last term in the first set of parentheses is missing.)
Now
\begin{align*}
& (e_{n + k} e_{n + k - 1} \cdots e_{k + 2})
      (e_{n + k + 1} e_{n + k} \cdots e_{k + 2})        \\
& \hspace*{3em} {\mbox{}}
  = (e_{n + k} e_{n + k - 1} \cdots e_{k + 3})
      (e_{n + k + 1} e_{n + k} \cdots e_{k + 4})
      (e_{k + 2} e_{k + 3} e_{k + 2})        \\
& \hspace*{3em} {\mbox{}}
  = z^{-1} (e_{n + k} e_{n + k - 1} \cdots e_{k + 3})
      (e_{n + k + 1} e_{n + k} \cdots e_{k + 4}) e_{k + 2}   \\
& \hspace*{3em} {\mbox{}}
  = z^{-1} (e_{n + k} e_{n + k - 1} \cdots e_{k + 3} e_{k + 2})
      (e_{n + k + 1} e_{n + k} \cdots e_{k + 4}).
\end{align*}
Next, we combine the second factor above with the third factor
in the original expression:
\begin{align*}
& (e_{n + k + 1} e_{n + k} \cdots e_{k + 4})
      (e_{n + k + 2} e_{n + k + 1} \cdots e_{k + 3})        \\
& \hspace*{3em} {\mbox{}}
  = (e_{n + k + 1} e_{n + k} \cdots e_{k + 5})
      (e_{n + k + 2} e_{n + k + 1} \cdots e_{k + 6})
      (e_{k + 4} e_{k + 5} e_{k + 4}) e_{k + 3}       \\
& \hspace*{3em} {\mbox{}}
  = z^{-1} (e_{n + k + 1} e_{n + k} \cdots e_{k + 5})
      (e_{n + k + 2} e_{n + k + 1} \cdots e_{k + 6})
      e_{k + 4} e_{k + 3}   \\
& \hspace*{3em} {\mbox{}}
  = z^{-1} (e_{n + k + 1} e_{n + k}
                 \cdots e_{k + 5} e_{k + 4} e_{k + 3})
      (e_{n + k + 2} e_{n + k + 1} \cdots e_{k + 6}).
\end{align*}
We proceed similarly through the rest of the factors.
After $l$ steps, the factor in position $l + 1$ is
$e_{n + k + l} e_{n + k + l - 1} \cdots e_{k + 2 l + 2}.$
After $n - 1$ steps, and using the definitions of $g_{n}^{(k)}$
and $f_{n}^{(k)},$
we get
\[
\sum_{j = 1}^r v_j g_n^{(k)} v_j^*
 = z^{- (n - 1)} g_{n - 1}^{(k + 1)}
\andeqn
\sum_{j = 1}^r v_j f_n^{(k)} v_j^* = f_{n - 1}^{(k + 1)}.
\]
Thus $f_{n - 1}^{(k + 1)} \in C^* \big( B_{n + k}, \, f_n^{(k)} \big),$
so, using the induction hypothesis,
\[
B_{2 n + k - 1} = C^* \big( B_{n + k}, \, f_{n - 1}^{(k + 1)} \big)
  \S C^* \big( B_{n + k}, \, f_n^{(k)} \big).
\]
In particular, $e_l \in C^* \big( B_{n + k}, \, f_n^{(k)} \big)$
for $l \leq 2 n + k - 2.$
It remains to show that
$e_{2 n + k - 1} \in C^* \big( B_{n + k}, \, f_n^{(k)} \big),$
which takes some work.

To start, we claim that, for any $m$ and $j,$
we have
\begin{align*}
& \big[ e_{j + 1} e_{j + 2} \cdots e_{m + j - 2} e_{m + j - 1} \big]
  g_m^{(j)}
  \big[ e_{m + j - 1} e_{m + j - 2} \cdots e_{j + 2} e_{j + 1} \big]
    \\
& \hspace*{18em} {\mbox{}}
 = z^{-2 (m - 1)} g_{m - 1}^{(j + 2)} e_{j + 1}.
\end{align*}
Let $b$ be the left hand side of this equation.
To prove the equation, start by using the relation~(\ref{Jones2})
a total of $m - 1$ times to show that the product of the first
factor above and the first factor in the expression for $g_m^{(j)}$
is $z^{- (m - 1)} e_{j + 1}.$
Then use the relation~(\ref{Jones3}) to commute $e_{j + 1}$ past
all but the last term
in the second factor in the expression for $g_m^{(j)},$
getting
\begin{align*}
b & = z^{- (m - 1)}
      \big( e_{m + j + 1} e_{m + j} \cdots e_{j + 3} \big)
       e_{j + 1} e_{j + 2}
       \big( e_{m + j + 2} e_{m + j + 1} \cdots e_{j + 3} \big)
    \\
& \hspace*{3em} {\mbox{}}
       \cdots
       \big( e_{2 m + j - 1} e_{2 m + j - 2} \cdots e_{m + j} \big)
       \big[ e_{m + j - 1} e_{m + j - 2} \cdots e_{j + 1} \big].
\end{align*}
For the same reason, one can commute $e_{j + 1} e_{j + 2}$ past all
but the last term in the next factor,
then $e_{j + 1} e_{j + 2} e_{j + 3}$ past all
but the last term in the factor after that, etc.
This gives
\begin{align*}
b & = z^{- (m - 1)}
     \big( e_{m + j + 1} e_{m + j} \cdots e_{j + 3} \big)
     \big( e_{m + j + 2} e_{m + j + 1} \cdots e_{j + 4} \big)
    \\
& \hspace*{1.5em} {\mbox{}}
     \cdots
     \big( e_{2 m + j - 1} e_{2 m + j - 2}
                         \cdots e_{m + j + 1} \big)
     \big[ e_{j + 1} e_{j + 2} \cdots e_{m + j} \big]
     \big[ e_{m + j - 1} e_{m + j - 2} \cdots e_{j + 1} \big].
\end{align*}
By repeated application of~(\ref{Jones2}),
the terms in square brackets combine to give
$z^{- (m - 1)} e_{j + 1},$
and the other terms are by definition $g_{m - 1}^{(j + 2)}.$
This completes the proof of the claim.

Next, we claim that if $m \geq 1$ and
a \ca\  $C \S B_N,$ for large $N,$
contains $B_{m + j}$ and $g_m^{(j)},$
then $C$ also contains $g_{m - 1}^{(j + 2)}.$
To see this, choose a
quasi-basis for $E_{j + 1} \colon B_{j + 1} \to B_j,$ say
$\big( (w_l, w_l^*) \big)_{1 \leq l \leq r},$
with $w_1, w_2, \ldots, w_r \in B_{j + 1} \S C.$
Again letting $b$ be the left hand side of the equation
in the previous claim,
it follows directly from the formula that $b \in C.$
Therefore also $\sum_{l = 1}^r w_l b w_l^* \in C.$
Also, $\sum_{l = 1}^r w_l e_{j + 1} w_l^* = 1$
(apply the defining equation of a quasi-basis to the element~$1$),
and $g_{m - 1}^{(j + 2)}$ commutes with each $w_l$
(because $g_{m - 1}^{(j + 2)}$ is a product of the \pj s
$e_{j + 3}, \, e_{j + 4}, \, \ldots, \, e_{2 m + j - 1}$).
Using the previous claim, we therefore get
\begin{align*}
\sum_{l = 1}^r w_l b w_l^*
& = z^{-2 (m - 1)}
    \sum_{l = 1}^r w_l g_{m - 1}^{(j + 2)} e_{j + 1} w_l^*   \\
& = z^{-2 (m - 1)}
     g_{m - 1}^{(j + 2)} \sum_{l = 1}^r w_l e_{j + 1} w_l^*
  = z^{-2 (m - 1)} g_{m - 1}^{(j + 2)}.
\end{align*}
Since $z^{-2 (m - 1)}$ is invertible, the claim follows.

Take $C = C^* \big( B_{n + k}, \, f_n^{(k)} \big).$
Recall that $g_n^{(k)}$ is a power of $z$ times $f_n^{(k)},$
and that we proved that $B_{2 n + k - 1} \S C.$
Apply the previous claim, first with
$m = n$ and $j = k,$ to get $g_{n - 1}^{(k + 2)} \in C$;
then with $m = n - 1$ and $j = k + 2,$
to get $g_{n - 2}^{(k + 4)} \in C$;
etc.;
finally,
with $m = 2$ and $j = 2 n + k - 4,$
to get $g_{1}^{(2 n + k - 2)} \in C.$
Since $g_{1}^{(2 n + k - 2)} = e_{2 n + k - 1},$
it follows that $B_{2 n + k} \S C.$
Thus $C^* \big( B_{n + k}, \, f_n^{(k)} \big) = B_{2 n + k}.$
This completes the induction step.

Taking $k = 0,$ we get $C^* ( B_{n}, \, f_n ) = B_{2 n},$
as desired.
This completes the proof.
\end{proof}

\begin{rmk}\label{R:Bimodules}
Let $1 \in A \subset B$ be an inclusion of unital \ca s
with a faithful conditional expectation $E \colon B \to A$
with index-finite type.
For right Hilbert $B$-modules $X$ and $Y,$
we denote by $\Hom_B (X, Y)$ the set of $B$-linear maps
$T$ from $X$ to $Y$ which are bounded in the sense
that
$\| \langle T \xi, T \xi \rangle_B\|
    \leq M \| \langle \xi, \xi \rangle_B \|$
for some positive
$M$ and which have adjoints $T^*$
with respect to the $B$-valued inner products.
We put $\End_B (X) = \Hom_B (X, X).$
We use similar notation ${}_A\Hom (X, Y)$ and
${}_A\End(X)$ for left Hilbert $A$-modules $X$ and $Y,$
${}_A\Hom_B (X, Y) = {}_A\Hom (X, Y) \cap \Hom_B (X, Y)$
for Hilbert $A$--$B$ bimodules, etc.
As discussed after Lemma~1.10 of~\cite{kw},
for a Hilbert $B$--$B$ bimodule $X,$ we can make the following
identifications, in which $A'$ is interpreted as those
maps commuting with the left action of $A,$ and similarly for $B'$:
\[
A' \cap \End_{B} (X) = {}_A\End_{B} (X)
\andeqn
B' \cap \End_{B} (X) = {}_B\End_{B} (X).
\]

We now interpret $B_l$ as a Hilbert module over $A$
using the composition $E_1 \circ E_2 \circ \cdots \circ E_l$
of conditional expectations,
and similarly over $B$ omitting the term $E_1.$
(See the example after Lemma~1.10 of~\cite{kw}.)
The basic construction $B_2 = C^* (B, e_1)$
for an inclusion $A \subset B$ of index-finite type
can then be identified as $\End_A (B),$
with the Jones \pj\  $e_1 = e_A$ being the conditional expectation $E,$
regarded as an operator on $B.$
(See Definition~2.1.2, Proposition~1.3.3,
and  Definition~2.2.10 of~\cite{wata}.)

If $\Index (E) \in A,$
Proposition~\ref{tower} then identifies
$\End_A (B_l)$ with $B_{2 l}$ and
$\End_B (B_{l + 1})$ with $B_{2 l + 1}$
as Hilbert $B$--$B$ bimodules.
This gives
\[
A' \cap B_{2 l} = {}_A\End_A (B_l),
\,\,\,\,\,\,
A' \cap B_{2 l + 1} = {}_A\End_B (B_{l + 1}),
\]
\[
B' \cap B_{2 l} = {}_B\End_A (B_l),
\andeqn
B' \cap B_{2 l + 2} = {}_B\End_A (B_{l + 1})
\]
for all $l.$
Note that these are all finite dimensional.
(Use Propositions 2.3.4 and 2.7.3 of~\cite{wata}.)
\end{rmk}

\begin{lem}\label{L:Gen}
Let $1 \in A \subset B$ be an inclusion of unital \ca s
with a conditional expectation $E \colon B \to A$
with index-finite type.
Then the basic construction $C^* (B, e)$ is equal to $B e B.$
\end{lem}

\begin{proof}
Let $\big( (u_j, u_j^*) \big)_{1 \leq j \leq n}$
be a quasi-basis for $E.$
Let $y \in C^* (B, e).$
Identifying $C^* (B, e)$ with $\End_A (B)$
as in Remark~\ref{R:Bimodules},
$y$ corresponds to a right $A$-module endomorphism $T$ of $B.$
Then for $x \in B$ we have
$T (x) = \sum_{j = 1}^n T (u_j E (u_j^* x)).$
Since $T (u_j) \in B,$ this formula shows that
$y = \sum_{j = 1}^n T (u_j) e u_j^* \in B e B.$
\end{proof}

The following is a C*~version of Lemma~1.2 of~\cite{PP2}.

\begin{lem}\label{L:Mlt}
Let $1 \in A \subset B$ be an inclusion of unital \ca s
with a conditional expectation $E \colon B \to A$
with index-finite type.
Following Definition~\ref{D:FiniteDepth} but with different notation,
let $e$ be the Jones projection for this inclusion,
and let $C^* (B, e)$ be the basic construction.
Then for every $x \in C^* (B, e)$ there is a unique $b \in B$ such that
$b e = x e.$
Moreover, if $D \S A$ and $x \in D' \cap C^* (B, e),$
then $b \in D' \cap B.$
\end{lem}

\begin{proof}
We first claim that if $b, c \in B$ satisfy $b e = c e,$
then $b = c.$
It suffices to consider the case $c = 0.$
If $b e = 0,$ then
$E (b^* b) e = e b^* b e = 0.$
Identifying $C^* (B, e)$ with $\End_A (B)$
as in Remark~\ref{R:Bimodules},
we see that in fact $E (b^* b) = 0,$
which implies $b = 0$ because $E$ is faithful.
This proves the claim.

Uniqueness in the statement of the lemma follows.

Let $F \colon C^* (B, e) \to B$ be the dual conditional expectation
(again, following Definition~\ref{D:FiniteDepth}).
We claim that if $x \in C^* (B, e)$ then
$b = \Index (E) F (x e)$ satisfies the conclusion of the lemma.
By Lemma~\ref{L:Gen},
it suffices to consider $x = y_1 e y_2$ with $y_1, y_2 \in B.$
Then $x e = y_1 E (y_2) e,$
so the formula of Proposition~2.3.2 of~\cite{wata} gives
$F (x e) = \Index (E)^{-1} y_1 E (y_2).$
Therefore
\[
b e = \Index (E) \big[ \Index (E)^{-1} y_1 E (y_2) \big] e
    = x e.
\]
This proves the claim.

It remains only to show that if $x \in D' \cap C^* (B, e),$
then $b \in D' \cap B.$
Let $a \in D.$
Then $a \in A,$ so $a$ commutes with $e$ and $x,$
and we get
\[
a b e = a x e = x e a = b e a = b a e.
\]
The first claim in the proof now implies that $a b = b a.$
This completes the proof.
\end{proof}

\begin{cor}\label{C:Ideal}
Let $1 \in A \subset B$ be an inclusion of unital \ca s
with a conditional expectation $E \colon B \to A$
with index-finite type.
Adopt the notation of Definition~\ref{D:FiniteDepth}.
Then for $k \geq 1$ we have
\[
(A' \cap B_k) e_k (A' \cap B_k)
 = (A' \cap B_{k + 1}) e_k (A' \cap B_{k + 1}).
\]
\end{cor}

\begin{proof}
Fix $k.$
Proposition~2.3.4 of~\cite{wata} implies that
$E_k \colon B_k \to B_{k - 1}$ has index-finite type.
We may therefore apply Lemma~\ref{L:Mlt} to the
inclusion of $B_{k - 1}$ in $B_k,$
getting $(A' \cap B_k) e_k = (A' \cap B_{k + 1}) e_k.$
Taking adjoints and applying the same lemma again,
we obtain the conclusion.
\end{proof}

If $A$ is any unital \ca,
we let $Z (A)$ denote the center of $A.$

\begin{lem}\label{L:DimZ}
Let $1 \in A \subset B$ be an inclusion of unital \ca s
with a conditional expectation $E \colon B \to A$
with index-finite type.
Assume that $Z (A)$ is finite dimensional.
Adopt the notation of Definition~\ref{D:FiniteDepth}.
Then for $k \geq 0$ we have
\[
\dim_{\C} (Z (A' \cap B_k)) \leq \dim_{\C} (Z (A' \cap B_{k + 2})),
\]
with equality \ifo\  $e_{k + 1}$ is full in $A' \cap B_{k + 2}.$
\end{lem}

\begin{proof}
Using Propositions 2.3.4 and~2.7.3 of~\cite{wata}, and induction,
we see that $B_m' \cap B_n$ is finite dimensional whenever
$0 \leq m \leq n.$
Since $e_{k + 1} [A' \cap B_{k + 2}] e_{k + 1}$ is a corner
in $A' \cap B_{k + 2},$
it therefore suffices to prove that
\[
\dim_{\C} (Z (A' \cap B_k))
  = \dim_{\C} (Z (e_{k + 1} [A' \cap B_{k + 2}] e_{k + 1} )).
\]
Since $e_{k + 1}$ commutes with everything in $B_k,$
the set $(A' \cap B_k) e_{k + 1}$ is a \ca,
and Lemma~\ref{L:Mlt} further implies that
$x \mapsto x e_{k + 1}$ is an isomorphism
$A' \cap B_k \to (A' \cap B_k) e_{k + 1}.$
Using the commutation relation at the first step,
\[
(A' \cap B_k) e_{k + 1}
= E_k (A' \cap B_{k + 1}) e_{k + 1}
= e_{k + 1} (A' \cap B_{k + 1}) e_{k + 1}
\]
at the second step,
Corollary~\ref{C:Ideal} at the third step,
and $e_{k + 1} \in A' \cap B_{k + 2}$ at the fourth step,
we further get
\begin{align*}
(A' \cap B_k) e_{k + 1}
& = e_{k + 1} (A' \cap B_k) e_{k + 1} (A' \cap B_k)
              e_{k + 1}    \\
& = e_{k + 1} \big[ (A' \cap B_{k + 1}) e_{k + 1}
      (A' \cap B_{k + 1}) \big] e_{k + 1}   \\
& = e_{k + 1} \big[ (A' \cap B_{k + 2}) e_{k + 1}
      (A' \cap B_{k + 2}) \big] e_{k + 1}
  = e_{k + 1} (A' \cap B_{k + 2}) e_{k + 1}.
\end{align*}
Thus $A' \cap B_k \cong e_{k + 1} (A' \cap B_{k + 2}) e_{k + 1}.$
In particular, their centers have the same dimension.
\end{proof}

The following proposition is a C*~analog of parts of
Theorem~4.6.3 of~\cite{GHJ}.

\begin{prp}\label{depth}
Let $1 \in A \subset B$ be an inclusion of unital \ca s,
and let $E \colon B \to A$ be a conditional expectation
with index-finite type such that $\Index (E) \in A.$
Assume that $Z (A)$ is finite dimensional.
Let $B_0 = A, \, B_1 = B, \, B_3, \, B_4, \, \ldots$
and $e_1, e_2, \ldots$ be as in Definition~\ref{D:FiniteDepth}.
Then the following conditions are equivalent:
\begin{enumerate}
\item\label{P_depth_1}
There is $k_0 \in \N$ such that
$(A' \cap B_{k_0}) e_{k_0} (A' \cap B_{k_0}) = A' \cap B_{k_0 + 1}.$
(That is, $E$ has finite depth
in the sense of Definition~\ref{D:FiniteDepth}.)
\item\label{P_depth_1b}
There is $k_0 \in \N$ such that
$(A' \cap B_{k}) e_{k} (A' \cap B_{k}) = A' \cap B_{k + 1}$
for all $k \geq k_0.$
\item\label{P_depth_2}
There is $L \in \N$ such that
$\dim_{\C} (Z (A' \cap B_{2 k})) \leq L$ for $k \geq 1.$
\item\label{P_depth_3}
There is $M \in \N$ such that 
$\dim_{\C} (Z (A' \cap B_{2 k + 1})) \leq M$ for $k \geq 1.$
\item\label{P_depth_4}
There is $L \in \N$ such that, for any $k \geq 1,$
the number of equivalence classes
of irreducible
$A$--$A$ Hilbert bimodules which appear in $B_{k}$ is at most $L.$
\item\label{P_depth_5}
There is $M \in \N$ such that, for any $k \geq 1,$
the number of equivalence classes of irreducible
$A$--$B$ Hilbert bimodules which appear in $B_{k}$ is  at most $M.$
\end{enumerate}
Moreover, the number $k_0$ in~(\ref{P_depth_1b})
can be chosen to be the depth of the inclusion,
as in Definition~\ref{D:FiniteDepth}.
\end{prp}

\begin{proof}
Assume~(\ref{P_depth_1}), for some $k_0.$
Using~(\ref{P_depth_1}) at the second step,
$e_{k_0} e_{k_0 + 1} e_{k_0} = \Index (E)^{-1} e_{k_0}$
(with $\Index (E)^{-1} \in A' \cap B_{k_0 + 1}$ and invertible)
at the third step,
and Corollary~\ref{C:Ideal} at the fifth step,
we get
\begin{align*}
1 & \in A' \cap B_{k_0 + 1}
    = (A' \cap B_{k_0 + 1}) e_{k_0} (A' \cap B_{k_0 + 1}) \\
  & = (A' \cap B_{k_0 + 1}) e_{k_0} e_{k_0 + 1} e_{k_0}
                              (A' \cap B_{k_0 + 1}) \\
  & \subset (A' \cap B_{k_0 + 1}) e_{k_0 + 1} (A' \cap B_{k_0 + 1})
    = (A' \cap B_{k_0 + 2}) e_{k_0 + 1} (A' \cap B_{k_0 + 2}).
\end{align*}
Since the last expression is an ideal in $A' \cap B_{k_0 + 2},$
we see by looking at the second last expression that
\[
(A' \cap B_{k_0 + 1}) e_{k_0 + 1} (A' \cap B_{k_0 + 1})
  = A' \cap B_{k_0 + 2}.
\]
Repeating this, we get
\[
(A' \cap B_k) e_k (A' \cap B_k) =  A' \cap B_{k + 1}
\]
for all $k \geq k_0.$
This is~(\ref{P_depth_1b}), and also proves the last statement.

Assuming~(\ref{P_depth_1b}),
Corollary~\ref{C:Ideal} and Lemma~\ref{L:DimZ} imply that
\[
\dim_{\C} (Z (A' \cap B_{k + 2}))
       = \dim_{\C} (Z (A' \cap B_k))
\]
for $k \geq k_0.$
This gives both (\ref{P_depth_2}) and~(\ref{P_depth_3}).

We now prove that (\ref{P_depth_2}) implies~(\ref{P_depth_1}).
The proof that (\ref{P_depth_3}) implies~(\ref{P_depth_1}) is
the same, and is omitted.
So suppose that $L = \sup_{k \geq 1} \dim_{\C} (Z (A' \cap B_{2 k}))$
is finite.
Lemma~\ref{L:DimZ} implies that
$k \mapsto \dim_{\C} (Z (A' \cap B_{2 k}))$ is nondecreasing.
So there exists $k$ such that
$\dim_{\C} (Z (A' \cap B_{2 k})) = \dim_{\C} (Z (A' \cap B_{2 k + 2})).$
Now Lemma~\ref{L:DimZ} implies that
$(A' \cap B_{2 k + 1}) e_{2 k + 1} (A' \cap B_{2 k + 1})
       = A' \cap B_{2 k + 1}.$
This proves~(\ref{P_depth_1}) with $k_0 = 2 k + 1.$

The equivalence of~(\ref{P_depth_2}) and~(\ref{P_depth_4})
follows by using Remark~\ref{R:Bimodules} to see that
$\dim_{\C} (Z (A' \cap B_{2 l}))$
and $\dim_{\C} (Z (A' \cap B_{2 l - 2}))$ are
the numbers of equivalence classes of irreducible
Hilbert $A$--$A$ bimodules which appear in
$B_l$ and $B_{l - 1}$ respectively.
The equivalence of~(\ref{P_depth_3}) and~(\ref{P_depth_5})
is similar.
\end{proof}

\begin{prp}\label{relative}
Let $1 \in A \subset B$ be an inclusion of unital \ca s,
and let $E \colon B \to A$ be a conditional expectation
with index-finite type and such that
$\Index (E) \in A.$
Assume that $Z (A)$ is finite dimensional.
Let $A \subset B \subset B_2$ be the basic
construction.
If $A \subset B$ has finite depth, then so does $B \subset B_2.$
\end{prp}

\begin{proof}
We will prove that
$A' \cap B_{2 l + 1} \cong B' \cap B_{2 l + 2}$ for $l \geq 0.$
By Remark~\ref{R:Bimodules},
we have $A' \cap B_{2 l + 1} = {}_A\End_B (B_{l + 1})$
and $B' \cap B_{2 l + 2} = {}_B\End_A (B_{l + 1}).$
Define a map $\Phi$ from ${}_A\End_B (B_{l + 1})$
to ${}_B\End_A (B_{l + 1})$ by
$\Phi(T)(x) = (T(x^*))^*$
for $T \in {}_A\End_B (B_{l + 1})$ and $x \in B_{l + 1}.$
It is obvious that $\Phi$ is an isomorphism from $A' \cap B_{2 l + 1}$
to $B' \cap B_{2 l + 2}.$
So Condition~(\ref{P_depth_2}) of Proposition~\ref{depth}
for the inclusion $A \subset B$
implies Condition~(\ref{P_depth_3}) of Proposition~\ref{depth}
for the inclusion $B \subset B_2.$
\end{proof}

\begin{prp}\label{R:Sector}
Let $1 \in A \subset B$ be an inclusion of simple unital C*-algebras,
with a conditional expectation $E \colon B \to A$
with index-finite type.
The inclusion $A \subset B$ has finite depth
in the sense of Definition~\ref{D:FiniteDepth}
\ifo\  it  has finite depth
in the sense of Definition~4.5 of~\cite{Iz}.
\end{prp}

\begin{proof}
We use the notation of Definition~\ref{D:FiniteDepth}.

It follows from Proposition 2.10.11 of~\cite{wata}
that Condition~(\ref{P_depth_2}) of Proposition~\ref{depth}
does not depend on the choice of the conditional expectation,
as long as it has index-finite type.
By Lemma~\ref{MinimalConditional}
we may therefore assume that $\Index (E)$ is a scalar,
so that Proposition~\ref{depth} applies.

Let $\io \colon A \otimes \cK \to B \otimes \cK$ be the inclusion.
For any sector $\et \in {\mathrm{Sect}} (B, A),$
in the sense of Section~4 of~\cite{Iz},
let $N (\et)$ denote the number of distinct irreducible sectors
in the decomposition of $\et$ of Lemma~4.1 of~\cite{Iz}.
We prove the following:
\[
N ([ (\io {\overline{\io}} )^k ])
  = \dim_{\C} (Z (B' \cap B_{2 k + 1})),
\,\,\,\,\,\,
N ([ (\io {\overline{\io}} )^k \io ])
  = \dim_{\C} (Z (A' \cap B_{2 k + 1})),
\]
\[
N ([ ({\overline{\io}} \io )^k ])
  = \dim_{\C} (Z (A' \cap B_{2 k})),
\andeqn
N ([ ({\overline{\io}} \io )^k {\overline{\io}}])
  = \dim_{\C} (Z (B' \cap B_{2 k + 2})).
\]
The result
will then follow from
Proposition~\ref{depth} and Proposition~\ref{relative}.

First, 
let $\gm \colon B \otimes \cK \to B \otimes \cK$
be the canonical endomorphism of Lemma~4.2 and Remark~4.3 of~\cite{Iz}.
Then Remark~4.3 of~\cite{Iz}
gives the downward tower of basic constructions
\[
(B \otimes \cK) \supset (A \otimes \cK) \supset \gm (B \otimes \cK)
 \supset \gm (A \otimes \cK) \supset \gm^2 (B \otimes \cK)
 \supset \gm^2 (A \otimes \cK) \supset \cdots.
\]
Reverse the segment of this tower
ending at $\gm^{k + 1} (A \otimes \cK).$
Using $\gm = \io {\overline{\io}},$
we can identify the result as
\[
A \otimes \cK \hookrightarrow B \otimes \cK
 \hookrightarrow A \otimes \cK
 \hookrightarrow \cdots \hookrightarrow A \otimes \cK
 \hookrightarrow B \otimes \cK,
\]
with all inclusions $A \otimes \cK \hookrightarrow B \otimes \cK$
given by $\io$ and all inclusions
$B \otimes \cK \hookrightarrow A \otimes \cK$
given by ${\overline{\io}}.$
Since the basic construction is preserved by tensoring with $\cK,$
this tower is isomorphic to
\[
A \otimes \cK \S B \otimes \cK \S B_2 \otimes \cK
 \S \cdots \S B_{2 k + 1} \otimes \cK \S B_{2 k + 2} \otimes \cK.
\]
So, for example,
\[
({\overline{\io}} \io )^k {\overline{\io}} (B \otimes \cK)' \cap
        M (A \otimes \cK)
 \cong (B \otimes \cK)' \cap M (B_{2 k + 2} \otimes \cK).
\]
Remark~2.12 of~\cite{Iz} shows that the right hand side is
isomorphic to $B' \cap B_{2 k + 2}.$
It is clear from the proof of Lemma~4.1 of~\cite{Iz}
that if $C$ and $D$ are any simple stable $\sm$-unital \ca s,
and if $\rh \colon C \to D$ defines a sector,
then $N ([\rh]) = \dim_{\C} (\rh (C)' \cap M (D)).$
The relation
$N ([ ({\overline{\io}} \io )^k {\overline{\io}}])
  = \dim_{\C} (Z (B' \cap B_{2 k + 2}))$
now follows.
The other three relations to be proved follow similarly.
\end{proof}

We will need two results about the cut-down $p \in p A p \subset p B p$
of an inclusion $1 \in A \subset B$ by
a \pj\  $p \in A' \cap B.$
(Note that $p A p$ is usually not contained in $A,$
because $p$ need not be in $A.$)
These are proved in the setting of factors in
Remark~2.6 of~\cite{Bi}.

The next proposition looks very much like Corollary~4.2 of~\cite{OT},
but differs in that
here the \pj\  $p$ is in $A' \cap B$ rather than $A.$

\begin{prp}\label{P:CutIFType}
Let $1 \in A \subset B$ be an inclusion of unital \ca s
with index-finite type,
and suppose $A$ is simple.
If $p \in A' \cap B$ is a nonzero projection, then
the inclusion $p A p \subset p B p$
also has index-finite type.
\end{prp}

\begin{proof}
Following Lemma~\ref{MinimalConditional},
let $E \colon B \to A$ be a pseudominimal conditional expectation.
Then $E (p)$ is a nonzero element in $A \cap A' = {\mathbb{C}},$
so $E (p) = \ld \cdot 1$ for some $\ld \in (0, \infty).$
Define a map $F$ from $p B p$ onto $p A p$ by
$F (x) = \ld^{-1} E (x) p$ for $x \in p B p.$
It is easy to see that $F$ is a conditional expectation
from $p B p$ onto $p A p.$
Let $\big( (u_j, u_j^*) \big)_{1 \leq j \leq n}$
be a quasi-basis for $E.$
We claim that
$\big( \big( \ld^{1/2} p u_j p, \,
            \ld^{1/2} p u_j^* p \big) \big)_{1 \leq j \leq n}$
is a quasi-basis for $F.$
For any $x \in B,$
we have $E (p x) = E (x p)$ by Definition~\ref{D:MinExpt}.
Using this and $p \in A'$ at the third step,
for any element $x \in p B p \subset B$ we have
\[
x = p x p
  = \sum_{j = 1}^n p u_j E (u_j^* p x p) p
  = \sum_{j = 1}^n p u_j p E (p u_j^* p x) p
  = \sum_{j = 1}^n
         \ld^{1/2} p u_j p F \big( \ld^{1/2} p u_j^* p x \big).
\]
A similar argument proves that
\[
x = \sum_{j = 1}^n
        F \big( x \ld^{1/2} p u_j p \big) \ld^{1/2} p u_j^* p.
\]
This proves the claim.
The existence of a quasi-basis implies that $F$ has index-finite type.
\end{proof}

\begin{prp}\label{P:CutFDepth}
Let $1 \in A \subset B$ be an inclusion of unital \ca s
with index-finite type and finite depth,
and suppose $A$ is simple.
If $p$ is a nonzero projection in $A' \cap B,$ then
the inclusion $p A p \subset p B p$
also has finite depth.
\end{prp}

\begin{proof}
Lemma~\ref{MinimalConditional} provides a conditional expectation
$E \colon B \to A$ with index-finite type
and such that $\Index (E)$ is a scalar.
Let
\[
A \subset B \subset B_2 \subset \cdots \subset B_k \subset \cdots
\andeqn
p A p \subset p B p
    \subset C_2 \subset \cdots \subset C_k \subset \cdots
\]
be the towers of iterated basic constructions.

Let $R_p$ be the projection in $\mathrm{End}_A (B) = B_2$
(see Remark~\ref{R:Bimodules} for this identification)
defined by
$R_p (b) = b p$ for $b \in B.$
Then $R_p \in B' \cap B_2,$
because the embedding of $B$ in $\mathrm{End}_A (B)$
is via left multiplications.
In particular $p R_p = R_p p.$
Let $p_1$ be
the projection in $A' \cap B_2$
defined by $p_1 = p R_p,$ that is,
as an operator in $\mathrm{End}_A (B),$
we have $p_1 (b) = p b p$ for $b \in B.$
Define an isomorphism $\varphi$ from $p_1 B_2 p_1$ onto
$\End_{p A p} (p B p)$
by $\varphi(p_1 x p_1) = p_1 x p_1 |_{p B p}.$
Then for $x, b \in B,$ we have $(p_1 x p_1) (p b p) = p x p (p b p),$
so $\varphi (p_1 B p_1) = p B p$
and $\varphi (p_1 A p_1) =  p A p.$
Thus we have an isomorphism of inclusions
\[
\left( p A p \subset p B p \subset C_2 \right)
\cong \left( p_1 A p_1 \subset p_1 B p_1 \subset p_1 B_2 p_1 \right).
\]

By Proposition~\ref{tower},
we can identify $\End_A (B_{2^{l - 1}})$ with $B_{2^l}$
for any $l \geq 2.$
The argument of the previous paragraph therefore gives \pj s
$p_l = p_{l - 1} R_{p_{l - 1}} \in A' \cap B_{2^l}$ such that
\[
p_l (x) = p_{l - 1} x p_{l - 1}
\]
for $x \in B_{2^{l - 1}},$
and isomorphisms of inclusions
\[
\left( p A p \subset C_{2^{l - 1}} \subset  C_{2^l} \right)
\cong
\left (p_l A p_l \subset p_l B_{2^{l - 1}} p_l
                \subset p_l B_{2^l} p_l \right).
\]
Since $A \subset B$ has finite depth,
Proposition~\ref{depth} provides
a constant $M$ such that
$\dim_{\C} ( Z (A' \cap B_l)) \leq M$ for $l \geq 2.$
Therefore we have
\[
\dim_{\C} \big( Z ((p A p)' \cap C_{2^l}) \big)
     = \dim_{\C} \big( Z ((p_l A p_l)' \cap p_l B_{2^l} p_l ) \big)
     \leq \dim_{\C} \big( Z (A' \cap B_{2^l}) \big)
     \leq M.
\]
Since
$\big(
 \dim_{\C} \big(Z ( (p A p)' \cap C_{2 n}) \big) \big)_{n \in \N}$
is nondecreasing (by Lemma~\ref{L:DimZ}),
we have
\[
\dim_{\C} (Z ((p A p)' \cap C_{2 n})) \leq M
\]
for $n \in \N.$
The algebra $p A p$ is simple (being isomorphic to $A$),
so Proposition~\ref{P:CutIFType} and Lemma~\ref{MinimalConditional}
provide a conditional expectation
$F \colon p B p \to p A p$ with index-finite type
and such that $\Index (F)$ is a scalar.
Now Proposition~\ref{depth} implies that
$p A p \subset p B p$ has finite depth.
\end{proof}

\section{Cancellation for inclusions}\label{S:Canc}

In this section we prove a cancellation theorem
for inclusions of simple \ca s with index-finite type.
We need the following modification of Blackadar's
cancellation theorem in~\cite{Bl7}.
(Also see Theorem~4.2.2 of~\cite{bl2}.)
Since that theorem is itself a modification
of an argument of Rieffel~\cite{Rf2},
we give a detailed proof for the reader's
convenience.
It is based on an argument of Goodearl~\cite{Kg},
which is given here with his permission.

\begin{thm}[Blackadar~\cite{Bl7}]\label{T:BCanc}
Let $A$ be a simple \ca.
Let $\cP \S \MIA$ be a set of \nzp s with the following two
properties:
\begin{enumerate}
\item\label{T:BCanc:1}
For every \nzp\  $q \in \MIA,$
there exists $p \in \cP$ such that $2 [p] \leq [q]$ in $K_0 (A).$
\item\label{T:BCanc:2}
$\sup_{p \in \cP} \tsr (p \MIA p) < \I.$
\end{enumerate}
Then the \pj s in $\MIA$ satisfy cancellation.
\end{thm}

In the arguments leading up to the proof,
we will tacitly identify \pj s
with their Murray-von Neumann equivalence classes.
For example, if $p, q \in \MIA$ are \pj s, we will let $p \oplus
q$ stand for any specific \pj\  in $\MIA$ with the appropriate
\mvnc\ class, without saying which one.
We further let $n * p$
stand for the direct sum of $n$ copies of $p.$
We use $\sim$ for \mvnc,
and we write $p \precsim q$ when there exists $r$ such that
$p \sim r \leq q.$

We will need the following result of Warfield,
which we restate here for \ca s and in terms of \pj s.

\begin{thm}\label{T:Warfield}
Let $A$ be a \ca,
let $p, q, r \in \MIA$ be \pj s, and let $n \in \N.$
Assume that
\[
p \oplus r \sim q \oplus r,
\,\,\,\,\,\,
\tsr (r \MIA r) \leq n,
\andeqn
n * r \precsim p.
\]
Then $p \sim q.$
\end{thm}

\begin{proof}
Using the Bass stable rank $\Bsr (A)$ in place of $\tsr (A),$
and using modules in place of \pj s,
this is Theorem~2.1 of~\cite{Rf2},
which is really just a combination
of Theorems~1.2 and~1.6 of~\cite{Wr}.
But $\Bsr (A) \leq \tsr (A)$ by Corollary~2.4 of~\cite{Rf1}.
(In fact, for any unital \ca~$A,$
we have $\Bsr (A) = \tsr (A)$ by~\cite{HV}.)
\end{proof}

We further need the following well known lemma,
whose proof we omit.

\begin{lem}\label{L:Dom}
Let $A$ be a simple \ca, and let $p, q \in \MIA$ be \pj s
with $q \neq 0.$
Then there exists $l \in \N$ such that $p \precsim l * q.$
\end{lem}


\begin{proof}[Proof of Theorem~\ref{T:BCanc}]
Set $n = \sup_{p \in \cP} \tsr (p \MIA p).$
By iterating Condition~(1) of the hypothesis,
we find that for every $m \in \N$ and
every \nzp\  $q \in \MIA,$
there exists $p \in \cP$ such that $2^m [p] \leq [q]$ in $K_0 (A).$
In particular,
for every \nzp\  $q \in \MIA,$
there exists $p \in \cP$ such that $n [p] \leq [q]$ in $K_0 (A).$

We now claim that
for every \nzp\  $q \in \MIA,$
there exists $p \in \cP$ such that $n * p \precsim q.$
To prove the claim,
given $q,$
choose $p \in \cP$ such that $n [p] \leq [q]$ in $K_0 (A).$
By the definition of the order on $K_0 (A),$
there are \pj s $r, s \in \MIA$ such that
$n * p \oplus r \oplus s \sim q \oplus s.$
Lemma~\ref{L:Dom} provides $l$ such that
$s \precsim l * p,$ whence
\begin{equation}\label{E:Sim}
n * p \oplus r \oplus l * p \sim q \oplus l * p.
\end{equation}

Let $l_0$ be the smallest nonnegative integer $l$
such that~(\ref{E:Sim}) holds.
The claim will be proved if we can show that $l_0 = 0.$
Suppose $l_0 > 0.$
Apply Theorem~\ref{T:Warfield}
with $(n + l_0 - 1) * p \oplus r$ in place of $p,$
with $q \oplus (l_0 - 1) * p$ in place of $q,$
and with $p$ in place of $r.$
We conclude that
\[
n * p \oplus r \oplus (l_0 - 1) * p \sim q \oplus (l_0 - 1) * p,
\]
contradicting the choice of $l_0.$
This completes the proof of the claim.

Now we can prove the theorem.
Let $p, q, r \in \MIA$ satisfy $p \oplus r \sim q \oplus r.$
Use the claim to choose $e \in \cP$
such that $n * e \precsim p.$
Thus, there is a \pj\  $s \in \MIA$ such that $n * e \oplus s \sim p.$
We then have $n * e \oplus s \oplus r \sim q \oplus r.$
Lemma~\ref{L:Dom} provides $l$ such that
$r \precsim l * e,$ whence
\begin{equation}\label{E:S2}
n * e \oplus s \oplus l * e \sim q \oplus l * e.
\end{equation}
As before, let $l_0$ be the least possible value of $l$
in~(\ref{E:S2}),
and if $l_0 > 0$ apply Theorem~\ref{T:Warfield}
with $(n + l_0 - 1) * e \oplus s$ in place of $p,$
with $q \oplus (l_0 - 1) * e$ in place of $q,$
and with $e$ in place of $r,$
getting
\[
n * e \oplus s \oplus (l_0 - 1) * e \sim q \oplus (l_0 - 1) * e.
\]
This is a contradiction, whence $l_0 = 0,$
and $p \sim n * e \oplus s \sim q.$
\end{proof}

We need several observations about the topological stable
rank of inclusions of \ca s to prove cancellation for crossed
products.
Recall that an inclusion $1 \in A \subset B$ of unital \ca s
is called {\emph{irreducible}} if
$A' \cap B = \C 1.$
(See Example~3.14 of~\cite{Iz}.)
The irreducible case of the following lemma
is essentially contained in Theorem~2.1 of~\cite{O}.

\begin{lem}\label{FirstCondition}
Let $1 \in A \subset B$ be an inclusion of
simple unital C*-algebras with index-finite type and finite depth.
Suppose that $A$ has property~(SP).
Then for any nonzero projection $q \in B$
there exists a nonzero projection $p \in A$
such that $p \precsim q.$
\end{lem}

\begin{proof}
The conclusion does not depend
on the choice of the conditional expectation,
so use Lemma~\ref{MinimalConditional} to choose $E \colon B \to A$
with index-finite type such that $\Index (E)$ is a scalar.

We first assume that $1 \in A \subset B$ is irreducible.
Let $q \in B$ be a nonzero projection.
Since the inclusion $1 \in A \subset B$ is irreducible
and has index-finite type and finite depth,
Proposition~\ref{R:Sector} and
Corollary~7.6 of~\cite{Iz} (also see Remark~2.4(3) of~\cite{Iz})
imply the following outerness condition on $E$:
for any $x \in B$ and for any nonzero hereditary C*-subalgebra C of $A,$
\[
\inf \big( \big\{ \| c (x - E (x)) c\|
    \colon {\mbox{$c \in C_{+},$ $\|c\| = 1$}} \big\} \big) = 0.
\]
This condition is equivalent to outerness in the sense
of Definition~2.2 of~\cite{O}.
By Theorem~2.1 and the proof of Corollary~2.3 of~\cite{O},
there is a nonzero projection $p \in A$
such that $p \precsim q.$

Now we consider the general case.
Since the inclusion has index-finite type, the relative commutant
$A' \cap B$ is finite dimensional.
(See Proposition 2.7.3 of~\cite{wata}.
The notation $C_B (A)$ there is the relative commutant,
defined in the statement of Proposition 1.2.9 of~\cite{wata}.)
Thus
\[
A' \cap B \cong \bigoplus_{l = 1}^n M_{k (l)} (\C)
\]
for suitable $n$ and $k (1), \, k (2), \, \ldots, \, k (n).$
For $1 \leq l \leq n,$
let $\big( e_{i, j}^{(l)} \big)_{1 \leq i, j \leq k (l)}$
be a system of matrix units
for the summand $M_{k(l)}.$
Each inclusion
$e_{j, j}^{(l)} A = e_{j, j}^{(l)} A e_{j, j}^{(l)} \subset
            e_{j, j}^{(l)} B e_{j, j}^{(l)}$
is irreducible,
and by Propositions \ref{P:CutIFType} and~\ref{P:CutFDepth}
has index-finite type and finite depth.
(Note: $e_{j, j}^{(l)} A e_{j, j}^{(l)}$ is usually not
contained in $A.$)
Set $k = \sum_{l = 1}^n k (l).$

Let $q \in B$ be a nonzero projection.
Theorem~5.1 of~\cite{O} implies that
$B$ has Property~(SP).
Since $B$ is simple,
by Lemma~3.2 of~\cite{LnTAF},
there is a projection $q_0 \in B$
such that, in the notation introduced before Theorem~\ref{T:Warfield},
we have $k * q_0 \precsim q.$
By Lemma~3.5.6 of~\cite{Hl},
there is a nonzero \pj\  $q_j^{(l)} \leq e_{j, j}^{(l)}$
such that $q_j^{(l)} \precsim q_0.$
The irreducible case provides
nonzero projections
\[
r_j^{(l)} \in e_{j, j}^{(l)} A e_{j, j}^{(l)}
\andeqn
s_j^{(l)} \in q_j^{(l)} e_{j, j}^{(l)} B e_{j, j}^{(l)} q_j^{(l)}
   = q_j^{(l)} B q_j^{(l)}
\]
such that $r_j^{(l)} \sim s_j^{(l)}.$
Thus $r_j^{(l)} \precsim q_0.$
Since $A$ is simple,
the map $A \to e_{j, j}^{(l)} A$ is bijective,
so there exists a \pj\  $p_j^{(l)} \in A,$
necessarily nonzero, such that
$r_j^{(l)} = e_{j, j}^{(l)} p_j^{(l)}.$

Since $A$ is simple and has Property~(SP),
repeated application of Lemma~1.8 of~\cite{PhtRp1a}
provides a nonzero projection
$p \in A$ such that $p \precsim p_j^{(l)}$
for all $l$ and $j.$
Then we have
\[
p = \sum_{l = 1}^n \sum_{j = 1}^{k (l)} p e_{j, j}^{(l)}
  \precsim
        \sum_{l = 1}^n \sum_{j = 1}^{k (l)} p_j^{(l)} e_{j, j}^{(l)}
  = \sum_{l = 1}^n \sum_{j = 1}^{k (l)} r_j^{(l)}
  \precsim k * q_0
  \precsim q.
\]
This completes the proof.
\end{proof}

\begin{prp}\label{P:Indexfinite}
Let $1 \in A \subset B$ be an inclusion
of simple unital \ca s with index-finite type and finite depth.
Suppose that $\tsr (A) = 1$ and $A$ has Property~(SP).
Then $B$ has cancellation.
\end{prp}

\begin{proof}
Let the notation be as in Definition~\ref{D:FiniteDepth},
and assume that $A \subset B$ has depth $m.$
Choose $k \in \N$ such that $k$ is odd and $k \geq m.$
Let $\cP$ be the set of all \nzp s in $A.$
We claim that, as a subset of $B_k,$ this set satisfies
the conditions of Theorem~\ref{T:BCanc}.

The \ca\  $B_2$ is stably isomorphic to $A.$
(In~\cite{wata}, see Proposition~1.3.4 and the argument
preceding Lemma~3.3.4.)
Continuing by induction,
and because $k + 1$ is even,
we find that $B_l$ is stably isomorphic to $A$ when $l$ is even,
and to $B$ when $l$ is odd.
Thus every $B_l$ is simple.
Moreover,
$B_{k + 1}$ is stably isomorphic to $A$ and
$B_{k}$ is stably isomorphic to $B.$
Theorem~6.4 of~\cite{Rf1} therefore gives $\tsr (B_{k + 1}) = 1.$
Also, by Proposition~1.6.6 of~\cite{wata} and induction,
every inclusion $B_l \subset B_{l + 1}$ has
index-finite type.

We claim that Condition~(\ref{T:BCanc:2}) of Theorem~\ref{T:BCanc}
is satisfied.
Proposition~\ref{depth} implies
\[
(A' \cap B_k) e_k (A' \cap B_k) = A' \cap B_{k + 1}.
\]
By Proposition~4.4 of~\cite{OT},
there are $n \in \N$ and $u_1, u_2, \ldots, u_n \in A' \cap B_{k + 1}$
such that, for every $p \in \cP,$
the family $\big( (p u_j, u_j^* p) \big)_{1 \leq j \leq n}$
is a quasi-basis for the conditional expectation
$F_p = E_k |_{p B_{k + 1} p}$ from $p B_{k + 1} p$ onto $p B_k p.$
For any such \pj~$p,$
using Corollary~\ref{Cor:Proposition5.3}
at the first step
and $\tsr (B_{k + 1}) = 1$ and Theorem~4.5 of~\cite{bl4}
at the second step,
we get
\[
\tsr (p B_k p)
 \leq n \cdot \tsr (p B_{k + 1} p) + n^2 - 2 n + 1
 = n^2 - n + 1.
\]
This proves the claim.

We next prove by induction on $l$ that $B_l$ has Property~(SP)
and that for every \nzp\  $q \in B_l$
there is a \nzp\  $p \in A$ such that $p \precsim q.$
The case $l = 0,$ corresponding to $B_l = A,$
is immediate.
So suppose the result is known for $l,$
and let $q \in B_{l + 1}$ be a \nzp.
The inclusion $B_l \subset B_{l + 1}$ has
finite depth by Proposition~\ref{relative}.
Both $B_l$ and $B_{l + 1}$ are simple,
so $B_{l + 1}$ has Property~(SP) by Theorem~5.1 of~\cite{O}.
Moreover,
Lemma~\ref{FirstCondition} then provides a \nzp\  $p_0 \in B_l$
such that $p_0 \precsim q,$
and the induction hypothesis provides a \nzp\  $p \in A$
such that $p \precsim p_0.$
This completes the induction.

To prove Condition~(\ref{T:BCanc:1}) for $\cP$
as a subset of $B_k,$
let $q \in M_{\infty} (B_k)$ be a \nzp.
Since $B_k$ has Property~(SP),
Lemma~3.5.6 of~\cite{Hl} provides
a \nzp\  $q_0 \in B_k$ such that
$q_0 \precsim q.$
The previous paragraph provides a \nzp\  $p_0 \in A$
such that $p_0 \precsim q_0.$
Since $A$ has Property~(SP),
by Lemma~3.5.7 of~\cite{Hl}
there exist orthogonal \mvnt\  \nzp s
$p_1, p_2 \in A$ with $p_1, p_2 \leq p_0.$
This completes the proof of Condition~(\ref{T:BCanc:1}).

We now conclude from Theorem~\ref{T:BCanc}
that $B_{k}$ has cancellation.
Since $B$ is stably isomorphic to $B_{k},$ so does $B.$
\end{proof}

Now we drop the requirement that the larger algebra be simple.
The following is our main theorem.

\begin{thm}\label{T:Mainindexfinite}
Let $1 \in A \subset B$ be an inclusion of unital \ca s of
index-finite type and with finite depth.
Suppose that $A$ is simple, $\tsr (A) = 1,$
and $A$ has Property~(SP).
Then $B$ has cancellation.
\end{thm}

\begin{proof}
Since $1 \in A \subset B$ has index-finite type,
results of~\cite{Iz} (Theorem~3.3, Definition~2.1,
and Remark~2.4(3) there)
provide projections $z_1, z_2, \ldots, z_k$ in the center of $B$
such that each $B z_j$ is simple and
\[
B = B z_1 \oplus B z_2 \oplus \cdots \oplus B z_k.
\]
By Propositions \ref{P:CutIFType} and~\ref{P:CutFDepth},
each inclusion $z_j \in z_j A z_j \subset B z_j$
has index-finite type and finite depth.
By Proposition~\ref{P:Indexfinite}, each $B z_j$ has cancellation.
Hence $B$ has cancellation.
\end{proof}

Using an observation by Blackadar and Handelman~\cite{bh}
we can sometimes determine $\tsr (B).$
Recall that a unital \ca~$A$ has {\emph{real rank zero}}
(see Theorem~2.6 of~\cite{bp0}) if every
selfadjoint element in $A$ can be approximated arbitrarily closely by
selfadjoint elements with finite spectrum.

\begin{cor}\label{C:StableIndexfinite}
Let $1 \in A \subset B$ be a pair of unital \ca s of
index-finite type and with finite depth.
Suppose that $A$ is simple with $\tsr (A) = 1$ and Property~(SP),
and that $B$ has real rank zero.
Then $\tsr (B) = 1.$
\end{cor}

\begin{proof}
The algebra $B$ has cancellation
by Theorem~\ref{T:Mainindexfinite}.
Since $B$ has real rank zero, $B$ has Property~(HP) by
Theorem~2.6 of~\cite{bp0}.
Therefore Theorem III.2.4 of~\cite{bh}
implies $\tsr (B) = 1.$
\end{proof}

\section{Cancellation for crossed products}\label{S:CrPrd}

In this section,
we apply the results of the previous section to crossed products by
finite groups,
in particular generalizing Theorem~5.4 of~\cite{jo}.
For comparison we
also give a result on crossed products by~${\mathbb{Z}}.$

The following result should be compared with
Question 8.2.3 of~\cite{bl3};
see Remark~\ref{rmk4.2b}.

\begin{cor}\label{C:Crossedproduct}
Let $A$ be an infinite dimensional
simple unital \ca, let $G$ be a finite group, and
let $\alpha$ be an action of $G$ on $A.$
Suppose that $\tsr (A) = 1$ and $A$ has Property~(SP).
Then $A \rtimes_{\af} G$ has cancellation.
Moreover, if $A \rtimes_{\af} G$ has real rank zero,
then $\tsr (A \rtimes_{\af} G) = 1.$
\end{cor}

\begin{proof}
Take $B = A \rtimes_{\af} G$ in
Theorem~\ref{T:Mainindexfinite}
and Corollary~\ref{C:StableIndexfinite}.
The finite depth condition is satisfied by Lemma~3.1 of~\cite{OT}.
\end{proof}

As an application of Corollary~\ref{C:Crossedproduct}
we get an interesting result when the
inclusion $A \subset B$ has index~$2,$
but does not necessarily have finite depth.

\begin{prp}\label{prp4.4}
Let $A$ be an infinite dimensional simple unital \ca\  %
with $\tsr (A) = 1$ and Property~(SP).
Suppose that the inclusion
$1 \in A \subset B$ has index $2.$
Then $B$ has cancellation.
\end{prp}

\begin{proof}
By Lemma 2.1.3 of~\cite{kt},
there is an action $\bt \colon {\mathbb{Z}} / 2 {\mathbb{Z}} \to B$
such that the basic construction $C^* (B, e_A)$ is isomorphic to
$B \rtimes_{\bt} {\mathbb{Z}} / 2 {\mathbb{Z}}.$
Proposition~1.3.4
and the discussion before Lemma 3.3.4 of~\cite{wata}
imply that $C^* (B, e_A)$ is stably isomorphic to $A.$
Therefore $C^* (B, e_A)$ is a simple unital
\ca\  with Property~(SP), and $\tsr (C^* (B, e_A)) = 1.$

Let
${\widehat{\bt}} \colon {\mathbb{Z}} / 2 {\mathbb{Z}}
   \to B \rtimes_{\bt} {\mathbb{Z}}/2{\mathbb{Z}}$
be the dual action.
By Takai duality (\cite{th}) we have
\[
\big( B \rtimes_{\bt} {\mathbb{Z}} / 2 {\mathbb{Z}} \big)
  \rtimes_{\widehat{\bt}} {\mathbb{Z}} / 2 {\mathbb{Z}}
\cong M_2 (B).
\]
Hence $B$ has cancellation by Corollary~\ref{C:Crossedproduct}.
\end{proof}

Let $\af \in \Aut (A)$
be an automorphism of a \ca~$A.$
There is no conditional expectation of index-finite
type from the crossed product $A \rtimes_{\af} {\mathbb{Z}}$ onto $A.$
Nevertheless, we have the following result.

\begin{thm}\label{thm4.4}
Let $A$ be a simple unital \ca\  with $\tsr (A) = 1$ and Property~(SP).
Let $\af \in \Aut (A)$ generate an outer action of ${\mathbb{Z}}$ on $A$
(that is, $\af^n$ is outer for every $n \neq 0$),
such that $\af_* = \id$ on $K_{0} (A).$
Then $A \rtimes_{\af} {\mathbb{Z}}$ has cancellation.
\end{thm}

\begin{proof}
Let $\cP$ be the set of all \nzp s in $\MIA,$
regarded as a subset of $M_{\infty} (A \rtimes_{\af} {\mathbb{Z}}).$
We claim that $\cP$ satisfies the conditions in Theorem~\ref{T:BCanc}.

For Condition~(\ref{T:BCanc:1}),
let $q \in M_{\infty} (A \rtimes_{\af} {\mathbb{Z}})$ be a \nzp.
Use Theorem~4.2 of~\cite{jo} to find
$e \in \cP$ such that $e \precsim q,$
and Proposition~3.5.6 of~\cite{Hl},
to find $p \in \cP$ such that $2[p] \leq [q].$

For Condition~(\ref{T:BCanc:2}),
let $p \in \cP.$
Assume $p \in M_n (A),$
and let $\af$ also denote the corresponding
automorphism of $M_n (A).$
Then $\alpha_{*} ([p]) = [p]$ because $\alpha_{*} = \id.$
Since
$\tsr (A) = 1,$ we can find a unitary $u \in M_n (A)$
such that $\af (p) = u^{*} p u.$
Define $\bt \in \Aut (M_n (A))$
by $\bt (a) = u \af (a) u^{*}.$
Then we calculate as follows,
in which the last two isomorphisms come from the proof of
Theorem~2.8.3(5) in \cite{Ph1}:
\[
p M_{\infty} (A \rtimes_{\af} {\mathbb{Z}}) p
 = p \big( M_n (A) \rtimes_{\af} {\mathbb{Z}} \big) p
 \cong p \big( M_n (A) \rtimes_{\bt} {\mathbb{Z}} \big) p
 \cong p M_n (A) p \rtimes_{\bt} {\mathbb{Z}}.
\]
Using Theorem~7.1 of~\cite{Rf1}, we therefore get
\[
\tsr \big( p M_{\infty} (A \rtimes_{\af} {\mathbb{Z}}) p \big)
  \leq \tsr (p M_n (A) p) + 1 = 2.
\]
Thus
$\sup_{p \in \cP}
      \tsr \big( p M_{\infty} (A \rtimes_{\af} {\mathbb{Z}}) p \big)
   \leq 2.$

Now Theorem~\ref{T:BCanc} implies that
$A \rtimes_{\af} {\mathbb{Z}}$ has cancellation.
\end{proof}

\begin{cor}\label{cor4.5}
Assume the hypotheses of Theorem~\ref{thm4.4}.
If in addition $A \rtimes_{\af} {\mathbb{Z}}$ has real rank zero,
then $\tsr (A \rtimes_{\af} {\mathbb{Z}}) = 1.$
\end{cor}

\begin{proof}
The proof is the same as that of
Corollary~\ref{C:StableIndexfinite}.
\end{proof}

\begin{exa}\label{exa4.1}
Let $A$ be a UHF~algebra or an irrational rotation algebra.
Let $\af \in \Aut (A)$ generate an outer action of ${\mathbb{Z}}.$
Then $A \rtimes_{\af} {\mathbb{Z}}$ has cancellation.
This follows from Theorem~\ref{thm4.4},
because $K_0 (A)$ has no nontrivial order preserving automorphisms.
\end{exa}

\begin{rmk}\label{rmk4.2a}
Let $A$ be a simple unital AT-algebra with real
rank zero and a unique tracial state.
Let $\af \in \Aut (A)$ be approximately inner
and generate an outer action of ${\mathbb{Z}}.$
Then the hypotheses of
Theorem~\ref{thm4.4} are satisfied.
In particular, by Corollary~\ref{cor4.5},
if $A \rtimes_{\af} {\mathbb{Z}}$ has real rank zero
then $\tsr (A \rtimes_{\af} {\mathbb{Z}}) = 1.$
In fact, under these hypotheses on $A$ and $\af,$
Kishimoto proved~\cite{kishi} that the following
are equivalent:
\begin{enumerate}
\item\label{rmk4.2a:1}
$A \rtimes_{\af} {\mathbb{Z}}$ has real rank zero.
\item\label{rmk4.2a:2}
$\af^m$ is uniformly outer for every $m \neq 0.$
\item\label{rmk4.2a:3}
$\af$ has the Rokhlin property.
\end{enumerate}
Theorem~\ref{thm4.4} suggests that, in this situation, one might deduce
$\tsr (A \rtimes_{\af} {\mathbb{Z}}) = 1$
{}from a weaker condition on $\af$ than the Rokhlin property.
We point out that Theorem~1.2 of~\cite{BEK} provides
many examples of actions of ${\mathbb{Z}}$ on the
$2^{\infty}$ UHF~algebra which do not have the Rokhlin
property and such that the crossed product is a simple AT~algebra
with real rank one (not zero), but of course with stable rank one.

\end{rmk}

\begin{rmk}\label{rmk4.2b}
Let $\af \colon G \to \Aut (A)$ be an action of a
discrete group $G$ on a unital \ca\  $A.$
Taking the crossed product $A \rtimes_{\af} G$ can
increase the topological stable rank if $G$ is finite and $A$ is
not simple (see Example 8.2.1 of~\cite{bl3}) or if $G$ is infinite
and $A$ is simple (see Example 8.2.2 of~\cite{bl3}).
Blackadar asked, in Question~8.2.3 of~\cite{bl3},
whether the crossed product of an AF~algebra by a
finite group has topological stable rank one.
This question remains open, even for simple AF~algebras and
${\mathbb{Z}} / 2 {\mathbb{Z}}.$
We have seen in Corollary~\ref{C:Crossedproduct}
that if $A$ is a simple unital \ca\  with
$\TR (A) = 0$ and $G$ is finite,
then $A \rtimes_{\af} G$ has cancellation.
It often happens that
cancellation for a simple unital \ca\  $B$ implies $\tsr (B) = 1,$
for example if $B$ has real rank zero.
However, a crossed product of a simple AF~algebra
by a finite group may have nonzero real rank (Example~9 of~\cite{el}),
and cancellation for a simple unital \ca\  $A$ does not imply
$\tsr (A) = 1$ (\cite{To}).
\end{rmk}

\end{document}